\documentclass[preprint,3p,12pt]{elsarticle}

\usepackage{color}
\usepackage{amssymb}
\usepackage{amsmath}
\usepackage{amsthm}

\usepackage{epsfig}
\usepackage{ulem}
\usepackage{epic}
\usepackage{tabls}
\usepackage{booktabs}
\usepackage{threeparttable}
\usepackage{pict2e}

\usepackage{float}
 \usepackage{subfig}

\usepackage{hyperref}

\hypersetup{backref,colorlinks=true,linkcolor=blue}

\usepackage[vlined]{algorithm2e}

\begin{document}
\begin{frontmatter}



\title{Multilevel  Second-Moment Methods
with  Group Decomposition \\ for Multigroup Transport Problems}

\author[ncsu,ncsu1]{Dmitriy Y. Anistratov}
\author[ncsu,ncsu2]{Joseph M. Coale}
 \author[lanl,lanl1]{James  S. Warsa}
 \author[lanl,lanl2]{ Jae H. Chang}
\address[ncsu]{Department of Nuclear Engineering,
North Carolina State University Raleigh, NC}
\address[ncsu1]{anistratov@ncsu.edu}
\address[ncsu2]{jmcoale@ncsu.edu}
\address[lanl]{Los Alamos National Laboratory,
Los Alamos, NM 87545, USA }
\address[lanl1]{warsa@lanl.gov}
\address[lanl2]{jhchang@lanl.gov}

\begin{abstract}
  This paper presents multilevel iterative schemes
   for solving the multigroup   Boltzmann transport  equations (BTEs)
    with parallel calculation of   group equations.
They are formulated with multigroup and grey low-order equations of the Second-Moment (SM) method.
The group high-order BTEs and low-order SM (LOSM) equations are solved in parallel.
To further improve convergence and increase computational efficiency of  algorithms
Anderson acceleration is applied to  inner iterations for solving the system of multigroup LOSM equations.
Numerical results are presented to demonstrate  performance of the  multilevel iterative methods.
 \end{abstract}


\end{frontmatter}

\section{Introduction}

The steady-state energy-dependent particle transport problems
are formulated by the multigroup Boltzmann transport equation (BTE) given by
\begin{equation} \label{bte}
 \boldsymbol{\Omega} \cdot \boldsymbol{\nabla} \psi_g(\boldsymbol{x},\boldsymbol{\Omega})   +
\sigma_{t,g}(\boldsymbol{x}) \psi_g(\boldsymbol{x},\boldsymbol{\Omega}) \! = \!
\frac{1}{4\pi}\sum_{g'=1}^{G} \! \sigma_{s,g'\to g}(\boldsymbol{x}
)\int_{4 \pi}\! \psi_{g'}(\boldsymbol{x},\boldsymbol{\Omega}') d \boldsymbol{\Omega}' \!+ \!
\frac{1}{4\pi} Q_g(\boldsymbol{x})   \, , \
g \in \mathbb{G} \, ,
\end{equation}
where   $\mathbb{G} = \{1,\ldots,G \}$.
Here standard notation is used.
This equation  models interaction of particles with matter  in a physical system with absorption and isotropic scattering.
It has  application for linear transport problems of various kind of particles,
for instance, neutrons, electrons,  and photons.
In nonlinear thermal radiative transfer (TRT) problems,
the time-dependent BTE is  coupled with the material energy balance (MEB) equation.
A group of methods for  TRT is based on linearization of the system of equations  on a time step.
This reduces the TRT problem to the BTE equation of the form \eqref{bte} with pseudo-scattering \cite{morel-jqsrt-1985,larsen-1988,morel-jcp-2007}.

Numerical transport algorithms for high performance computers use a variety of
approaches to achieve  efficient parallel computations for solving the linear BTE
\cite{hanus-ragusa-m&c2019,mla-jcp-2020}.
A natural  element  of parallel algorithms
is to perform calculations of group equations in parallel taking advantage of the particle transport problem's multigroup structure.
This can be interpreted as problem decomposition over one element of the phase space, namely, particle energy.
Efficient Diffusion Synthetic Acceleration (DSA) algorithms for multigroup transport problems can be formulated
with  decoupled group equations  \cite{jmc-jsw-dya-jhc-m&c2021}.

In this paper, we describe new  iterative methods for  multigroup  transport problems
 which  solve  the group equations  in parallel.
 They are formulated on the basis of
 low-order equations of the Second-Moment (SM) method
and nonlinear projective approach \cite{lewis-miller-1976,mla-ewl,dya-vyag-nse-2011}.
 The  low-order SM (LOSM) equations are similar to those of the  DSA  method \cite{mla-ewl}.
 The main difference is that the   LOSM   equations are formulated
 for the   angular moments of the solution,
 while the low-order DSA equations are defined for the iterative corrections of the moments.
Thus,  computational tools  based on the DSA  can be modified to use the SM method.
 We present a nonlinear multilevel SM (MLSM) method that
 consists of (i) multigroup  high-order  BTEs for group angular fluxes,
 (ii) multigroup LOSM equations for group scalar fluxes and currents,
 and (iii)  grey  LOSM equations for total scalar flux  and current.
 The scattering terms in both group high-order and LOSM equations are
 formulated in a nonlinear form.
The effective grey LOSM problem is defined  by means of     cross sections averaged
 with the  iterative group LOSM solution.
The group high-order BTEs and  group LOSM equations are solved in parallel at corresponding stages of  iteration algorithms.
The group-to-group   scattering terms in the group LOSM equations  are defined with lagged iterative solution.
In this case  both downscattering and upscattering of particles have a similar effect on
convergence of inner iterations with respect to energy groups.
Iterations of this kind can yield slow convergence without acceleration. The convergence rate of inner iterations over the system of multigroup low-order equations is improved with use of the grey LOSM equations.
To further accelerate convergence and increase computational efficiency of  parallel algorithms
we apply Anderson acceleration to the inner multigroup iterations \cite{anderson-1965}.

The reminder of the paper is organized as follows.
 In Sec. \ref{sec:mlsm}, the  MLSM  method is formulated.
 In Sec. \ref{sec:mlsm-aa}, we present the  MLSM  method  with Anderson Acceleration
 of inner iterations over multigroup  LOSM  equations.
 The numerical results are presented in Sec. \ref{sec:res}.
 We conclude with a discussion in Sec. \ref{sec:conc}.

\section{\label{sec:mlsm} Multilevel Second-Moment Method}

We consider transport problems in 1D slab geometry.
The  iteration scheme of the MLSM method  with groups solved in parallel is presented in Algorithm~\ref{alg-mlsm},
where
$\ell$ is the index of outer transport iterations,
 $k$ is the index of   the inner iterations between    multigroup  and grey   LOSM  equations, and
$s$ is the index of  the innermost  iterations for solving  multigroup  LOSM equations.
$k_{max}$ and $s_{max}$ are  the maximum numbers of the corresponding inner iterations.

\clearpage

\begin{algorithm}[t]
\DontPrintSemicolon
$\ell=-1$, $\psi_{g}^{(0)}=const$, $P_{g}^{(0)}=0$, $P^{(0)}=0$\;
\While{$|| \phi^{(\ell)}- \phi^{(\ell-1)} ||>\epsilon$ }{
 $\ell=\ell+1$ \;
 \If{$\ell>0$}{
 $\phi_g^{(\ell-1)}  \ \Rightarrow  \ \bar\sigma_{s,g}^{(\ell-1)}, \ g \in \mathbb{G}$  \;
 \For{all $g \in \mathbb{G}$ in parallel}{
    Level 1: Solve the high-order transport equation (Eqs. \eqref{ho-mlsm})  for group $g$
 $ \ \Rightarrow  \ \psi_{g}^{(\ell)}$\;}
 $\psi_{g}^{(\ell)} \ \Rightarrow
 \ P_{g}^{(\ell)} \  \mbox{for} \ g \in \mathbb{G}, \ P^{(\ell)}$\;
 }
   $k=0$\;
 \While{$k \le k_{max}$ }{
  $k=k+1$ \;
  $s=0$ \;
\While{$s \le s_{max}$ }{
 $s=s+1$ \;
 $\phi_{g}^{(s-1,k,\ell)}\, , \ \phi^{(k-1,\ell)}  \ \Rightarrow  \ \zeta^{(s-1,k-1,\ell)} $ \;
  \For{all $g \in \mathbb{G}$ in parallel}{
   Level 2: Solve the   LOSM equation  (Eqs.  \eqref{glosm-alg-1}) for group $g$
 $ \ \Rightarrow  \ \phi_{g}^{(s,k,\ell)}$, $J_{g}^{(s,k,\ell)}$\;
 }
 }
 $\phi_{g}^{(k,\ell)}  \leftarrow  \phi_{g}^{(s_{max},k,\ell)}$ , $J_{g}^{(k,\ell)}  \leftarrow J_{g}^{(s_{max},k,\ell)}$ \;
  $\phi_g^{(k,\ell)}$, $J_g^{(k,\ell)} \ \Rightarrow \ \bar\sigma_{a}^{(k,\ell)}$, $\bar \sigma_{t}^{(k,\ell)}$, $\bar \eta^{(k,\ell)}$ \;
Level 3: Solve the grey LOSM   equations  (Eqs. \eqref{grey-losm-it})
$\ \Rightarrow \ \phi^{(k,\ell)}$, $J^{(k,\ell)}$\;
}
$\phi_{g}^{(\ell)}  \leftarrow  \phi_{g}^{(k_{max},\ell)}$ , $J_{g}^{(\ell)}  \leftarrow  J_{g}^{(k_{max},\ell)}$,
$\phi^{(\ell)}  \leftarrow  \phi^{(k_{max},\ell)}$ , $J^{(\ell+1)}  \leftarrow  J^{(k_{max},\ell)}$ \;
}
\caption{\label{alg-mlsm}  The  MLSM  method  with  group equations solved in parallel}
\end{algorithm}

The multilevel hierarchy of equations of the MLSM method is defined as follows.
\begin{itemize}
\item{\bf Level 1.} The multigroup high-order transport equations  with decoupled groups
are given by
\begin{equation} \label{ho-mlsm}
\mu \frac{\partial  \psi_g^{(\ell)}}{\partial x} +
\sigma_{t,g} \psi_g^{(\ell)}=
\frac{1}{2} \bar \sigma_{s,g}^{(\ell-1)} \phi^{(\ell-1)} +
\frac{1}{2} Q_g   \, ,
\quad
\mbox{where}
\quad
\bar \sigma_{s,g}^{(\ell-1)}=
\frac{\sum_{g'=1}^{G}\sigma_{s,g'\to g}\phi_{g'}^{(\ell-1)}}{ \sum_{g'=1}^{G} \phi_{g'}^{(\ell-1)}} \, .
\end{equation}
The right-hand side (RHS) of the BTE  is transformed by means of (i) the total scalar flux that is the solution of the grey LOSM problem
and (ii) averaged cross section $\bar \sigma_{s,g}$  defined with the group scalar fluxes obtained from  the multigroup LOSM equations   \cite{dya-vyag-nse-2011}.

\item {\bf Level 2.}  The multigroup  LOSM  equations are defined by
  \begin{subequations}  \label{glosm-alg-1}
 \begin{equation}
 \frac{d    J_g^{(s,k,\ell)}}{d x}
+ \Big(\sigma_{t,g}- \sigma_{s,g\to g}\Big)    \phi_g^{(s,k,\ell)}=
   \zeta^{(s-1,k-1,\ell)}  \sum_{\stackrel{g'=1}{g'\ne g}}^{G}\sigma_{s,g'\to g}     \phi_{g'}^{(s-1,k,\ell)} +   Q_g \, ,
\end{equation}
 \begin{equation}
\frac{1}{3}  \frac{d    \phi_g^{(s,k,\ell)}}{d x}
 + \sigma_{t,g}    J_g^{(s,k,\ell)} =
 \frac{d P_g^{(\ell)}}{d x}   \, ,
\end{equation}
\end{subequations}
where
   \begin{equation}
  \zeta^{(s-1,k-1,\ell)}   = \frac{\phi^{(k-1,\ell)}}{\sum_{g'=1}^{G}   \phi_{g'}^{(s-1,k,\ell)}}\, ,
\quad
P_g^{(\ell)} = \int_{-1}^1 \Big( \frac{1}{3}   - \mu^2 \Big)  \psi_g^{(\ell)}  d \mu \, .
 \end{equation}
 The RHS   of Eq. \eqref{glosm-alg-1} is formulated
using the multiplicative correction factor $\zeta$
that is defined by the solution of the multigroup and grey LOSM equations  \cite{dya-vyag-nse-2011,dya-lrc-jpj-jcp-2017}.
This form of the RHS introduces nonlinearity in  the multigroup LOSM equations.

\item{\bf Level 3.} The grey LOSM equations have the following form:
\begin{subequations}\label{grey-losm-it}
 \begin{equation}
  \frac{d J^{(k,\ell)}}{d x}
+  \bar \sigma_{a}^{(k,\ell)}  \phi^{(k,\ell)}= Q \, ,
\end{equation}
 \begin{equation}
 \frac{1}{3}  \frac{d \phi^{(k,\ell)}}{d x}  +
\bar \sigma_t^{(k,\ell)} J^{(k,\ell)}  +
 \bar \eta^{(k,\ell)}\phi^{(k,\ell)}=
 \frac{d P^{(\ell)}}{d x} \, ,
\end{equation}
\end{subequations}
 where
\begin{equation}
\bar \sigma_{a}^{(k,\ell)} =
\frac{\sum_{g=1}^{G}\sigma_{a,g}\phi_{g}^{(k,\ell)}}{ \sum_{g=1}^{G} \phi_{g}^{(k,\ell)}} \, ,
\quad
 \bar \sigma_{t}^{(k,\ell)}  = \frac{\sum_{g=1}^{G}\sigma_{t,g} \Big|J_{g}^{(k,\ell)}\Big|}{ \sum_{g=1}^{G} \Big|J_{g}^{(k,\ell)}\Big|} \, ,
 \end{equation}
  \begin{equation}
\bar \eta^{((k,\ell)}  = \frac{\sum_{g=1}^{G}\Big(\sigma_{t,g} - \bar\sigma_{t}\Big)J_g^{(k,\ell)}}{ \sum_{g=1}^{G}   \phi_{g}^{(k,\ell)}} \, ,
\quad
P^{(\ell)} = \sum_{g=1}^{G} P_g^{(\ell)}  \, , \quad Q = \sum_{g=1}^{G} Q_g \, .
 \end{equation}
\end{itemize}
The high-order BTE (Eq. \eqref{ho-mlsm}) is discretized by the  linear-discontinuous (LD)  finite element method.
The spatial discretization of the multigroup and grey LOSM equations are
consistent with the LD transport scheme \cite{dya-jsw-nse-2018}.

\section{\label{sec:mlsm-aa}  The MLSM Method with Anderson Acceleration}

\subsection{Anderson Acceleration}

Let us  consider  a general  equation of the  following form:
\begin{equation} \label{G-eq}
\boldsymbol{\varphi}  = \mathcal{A}(\boldsymbol{\varphi})\, , \quad
\boldsymbol{\varphi} \in \mathbb{R}^n
\end{equation}
that is solved with the fixed-point iterations (FPI)
\begin{equation}
\boldsymbol{\varphi}^{(s+1)}  =  \mathcal{A} \big(\boldsymbol{\varphi}^{(s)}\big) \, .
\end{equation}
The residual for the $s^{th}$ iterate $\boldsymbol{\varphi}^{(s)}$ is defined by
\begin{equation}
\boldsymbol{r}(\boldsymbol{\varphi}^{(s)}) = \mathcal{A}(\boldsymbol{\varphi}^{(s)}) - \boldsymbol{\varphi}^{(s)} \, .
 \end{equation}
Anderson acceleration  is an iterative algorithm
that applies an extrapolation based on a linear combination of iterates and  the values of $\mathcal{A}$.
The coefficients of the linear combination are determined in such a way that they
minimize the linear combination of the corresponding sequence of residuals.
Algorithm \ref{alg-a}  presents the iteration scheme of Anderson acceleration   \cite{anderson-1965,walker-2011,toth-kelly-2015}.
The parameter $m$  defines  the maximum  algorithmic depth.
The set  of mixing  parameters $\beta_s$ are used for relaxation.
 We refer to this iteration algorithm as AA(m).
 \begin{algorithm}[ht]
\DontPrintSemicolon
Define  $\boldsymbol{\varphi}^{(0)}$ \;
$\boldsymbol{\varphi}^{(1)} = \mathcal{A}(\boldsymbol{\varphi}^{(0)}) \, , \
 \boldsymbol{r}(\boldsymbol{\varphi}^{(0)}) = \mathcal{A}(\boldsymbol{\varphi}^{(0)}) - \boldsymbol{\varphi}^{(0)}$\;
\For{s = 1,2,...}
{
$m_s = \min(m,s)$\;
 $\boldsymbol{r}(\boldsymbol{\varphi}^{(s)}) = \mathcal{A}(\boldsymbol{\varphi}^{(s)}) - \boldsymbol{\varphi}^{(s)}$\;
$
\min_{\boldsymbol{\alpha}^s} \big|\big| \sum_{j=0}^{m_s}\alpha_j^s \boldsymbol{r}(\boldsymbol{\varphi}^{(s-m_s +j)})\big|\big|_{2} \quad
\mbox{s. t.} \quad \sum_{j=0}^{m_s} \alpha_j^{s} = 1 \, .
$ \;
$
\boldsymbol{\varphi}^{(s+1)} = (1-\beta_s) \sum_{j=0}^{m_s} \alpha_j^{s} \boldsymbol{\varphi}^{(s-m_s+j)}
+  \beta_s  \sum_{j=0}^{m_s} \alpha_j^{s}  \mathcal{A}\big(\boldsymbol{\varphi}^{(s-m_s+j)}\big)
$
}
\caption{\label{alg-a} Anderson acceleration method for solving $\boldsymbol{\varphi}  = \mathcal{A}(\boldsymbol{\varphi})$}
\end{algorithm}

In this study, we use AA(1)  with    $\beta_s=1 $
that  converges r-linearly in $\ell_2$-norm  provided that the coefficients $\alpha_j^s$ are bounded \cite{toth-kelly-2015}.
This scheme
defines the next iterate as follows:
\begin{equation}
\boldsymbol{\varphi}^{(s+1)} = \alpha_0^{s}\mathcal{A}(\boldsymbol{\varphi}^{(s-1)}) + \alpha_1^{s}\mathcal{A}(\boldsymbol{\varphi}^{(s)})=
  \alpha_0^{s}\boldsymbol{\varphi}^{(s-1)} + \alpha_1^{s}\boldsymbol{\varphi}^{(s)}
+   \alpha_0^{s}\boldsymbol{r}(\boldsymbol{\varphi}^{(s-1)}) + \alpha_1^{s}\boldsymbol{r}(\boldsymbol{\varphi}^{(s)}) \, ,
\end{equation}
where $\boldsymbol{\alpha}^s=(\alpha_0^s , \alpha_1^s)^\top$ solves
\begin{equation} \label{min}
\min_{\boldsymbol{\alpha}^s} ||\alpha_0^s \boldsymbol{r}(\boldsymbol{\varphi}^{(s-1)}) +
\alpha_1^{s} \boldsymbol{r}(\boldsymbol{\varphi}^{(s)} ) \bigr)||_{2} \quad
\mbox{s. t.} \quad \alpha_0^{s} + \alpha_1^{s} = 1 \, .
\end{equation}
To determine   $\alpha_0^{s}$, we apply the following condition:
\begin{equation}
\Bigl( ||\boldsymbol{r}(\boldsymbol{\varphi}^{(s)})
 + \alpha_0^{s}  \bigl( \boldsymbol{r}(\boldsymbol{\varphi}^{(s-1)}) - \boldsymbol{r}\bigl(\boldsymbol{\varphi}^{(s)} ) \bigr)||_{2}^2 \Bigr)^{\prime}_{\alpha_0}= 0 \,.
 \end{equation}
This yields
\begin{equation}
 \alpha_0^{s}  =
 \frac{\sum_{i=1}^n r_i^{(s)}\bigl(r_i^{(s)} - r_i^{(s-1)}\bigr)}{
 \sum_{i=1}^n \bigl(r_i^{(s-1)} - r_i^{(s)}\bigr)^2} \, ,
\end{equation}
 where  $\{r_i^{(s)}\}_{i=1}^n = \{r_i(\boldsymbol{\varphi}^{(s)})\}_{i=1}^n $ .

\subsection{MLSM Algorithm with Anderson Acceleration of Innermost Iterations}

We apply  AA(1)  to the inner iterations of the multigroup LOSM equations at Level 2 (see Sec. \ref{sec:mlsm}).
The vector of the solution   $\boldsymbol{\varphi}=\{ \boldsymbol{\varphi}_g \}_{g=1}^G $
consists of $\boldsymbol{\varphi}_g$   defined by the grid functions of  $\phi_{g}$ and $J_g$. The residual is given by
$\boldsymbol{r}(\boldsymbol{\varphi})  = \{ \boldsymbol{r}_g\}_{g=1}^G $,
where $\boldsymbol{r}_g(\boldsymbol{\varphi} ) =\mathcal{F}_g \boldsymbol{\varphi}$ and $\mathcal{F}_g$
is the operator of the discretized LOSM equations in the group $g$. The  iteration scheme for the MLSM method
with AA(1) for the inner multigroup iterations is presented in Algorithm \ref{alg-MIM-v1c}.
Hereafter we refer to this algorithm as the MLSM-AA(1) method.

\begin{algorithm}[t]
\DontPrintSemicolon
 \ldots\;
$k=0$, \
  $\phi_{g}^{(0,\ell)}= \int_{-1}^{1}\psi_{g}^{(\ell)} d\mu$,
    $J_{g}^{(0,\ell)}=\int_{-1}^{1}\mu \psi_{g}^{(\ell)} d\mu $\;
  \While{$k \le k_{max}$ }{
  $k=k+1$, \   $s=0$ \;
  $\hat \phi_{g}^{(0,k,\ell)}= \phi_{g}^{(k-1,\ell)}$,
   $\hat J_{g}^{(0,k,\ell)}= J_{g}^{(k-1,\ell)}$,   $\phi_{g}^{(0,k,\ell)}= \phi_{g}^{(k-1,\ell)}$\;
   Calculate residual $\boldsymbol{r}(\boldsymbol{\hat \varphi}^{(0,k,\ell)})$
 of  the multigroup  LOSM equations\;
\While{$s \le s_{max}$ }{
 $s=s+1$ \;
  $\phi_{g}^{(s-1,k,\ell)}\, , \ \phi^{(k-1,\ell)}  \ \Rightarrow  \ \zeta^{(s-1,k-1,\ell)} $ \;
    \For{all $g \in \mathbb{G}$ in parallel}{
 Level 2: Solve the   LOSM equation   for group $g$
 $ \Rightarrow  \  \hat \phi_{g}^{(s,k,\ell)}$, $ \hat J_{g}^{(s,k,\ell)}$\;
 }
 Calculate residual $\boldsymbol{r}(\boldsymbol{\hat \varphi}^{(s,k,\ell+1)})$
 of  the multigroup  LOSM equations\;
$\boldsymbol{r}(\boldsymbol{\hat \varphi}^{(s,k,\ell)})\, , \ \boldsymbol{r}(\boldsymbol{\hat \varphi}^{(s-1,k,\ell)})  \ \Rightarrow \ \alpha_0 \, , \alpha_1$\;
$\boldsymbol{\varphi}^{(s,k,\ell)} = \alpha_0 \boldsymbol{\hat \varphi}^{(s-1,k,\ell)} +
\alpha_1 \boldsymbol{\hat \varphi}^{(s,k,\ell)}  +
\alpha_0 \mathbf{r}(\boldsymbol{\hat \varphi}^{(s-1,k,\ell)}) +
\alpha_1 \mathbf{r}(\boldsymbol{\hat \varphi}^{(s,k,\ell)})$ \;
$\boldsymbol{\varphi}^{(s,k,\ell)}  \ \Rightarrow \
\phi_{g}^{(s,k,\ell)}\, , \ J_{g}^{(s,k,\ell)}$\;
}
$\ldots$\;
}
$\ldots$\;
\caption{\label{alg-MIM-v1c}  The  MLSM-AA(1)  method    with  group equations solved in parallel.}
\end{algorithm}

\section{\label{sec:res} Numerical Results}

{\bf Test 1}.
This is a 10-group problem for a slab $0 \le x \le 32$ \cite{jmc-jsw-dya-jhc-m&c2021}. The cross sections are given in Table \ref{test1-xs} (see \ref{appendix}).
The groups are coupled with each other  due to downscattering and upscattering. The boundary conditions are vacuum. The external source is constant and $Q_g=1$ $\forall g$.
The spatial mesh is uniform with 128 cells.
There are  16  angular  directions. The double $S_8$ Gauss-Legendre  quadrature set is used.
The  convergence criterion is $\varepsilon =10^{-9}$.
Table \ref{test1-strength} shows the measure of connection strength of groups  given by \cite{mg-tutorial}
\begin{equation}
S_{gg'} = \frac{a_{gg'}}{\max_{g''\ne g} (a_{gg''})} \, , \quad
a_{gg'} = \sigma_{s,g '\to g} \, .
\end{equation}
In this  test, most of groups are strongly connected with other groups.
The group scattering is   high and in the range $0.9 \le c_g \le 0.9999$ (see Table \ref{test1-xs}).
The target number of   transport iterations in  Test 1 equals 15.

\begin{table}[t]
		\centering
\caption{\label{test1-strength} Connection strength of  groups ($S_{gg'}$) in Test 1}
{\small
\begin{tabular}{|c|c|c|c|c|c|c|c|c|c|c|}
\hline
$g \, \backslash  \, g'$	&	1	&	2	&	3	&	4	&	5	&	6	&	7	&	8	&	9	&	10	 \\  \hline
1	&	0	&	0	&	0	&	0	&	0	&	0	&	0	&	0	&	0	&	0	 \\
2	&	1. 	&	0 	&	0	&	0	&	0	&	0	&	0	&	0	&	0	&	0	 \\
3	&	0.71	&	1. 	&	0	&	0	&	0	&	0	&	0	&	0	&	0	&	0	 \\
4	&	1. 	&	0.53	&	0.36	&	0	&	0	&	0	&	0	&	0	&	0	&	0	 \\
5	&	1. 	&	0.21	&	0.48	&	1.	&	0	&	0	&	0	&	0	&	0	&	0	 \\
6	&	0	&	0.78	&	0.64	&	1. 	&	0.81	&	0	&	0	&	0	&	0	&	0	 \\
7	&	0	&	0	&	0.26	&	0.09	&	0.20	&	0.48	&	0	&	1.	&	0.22	&	0.23	 \\
8	&	0	&	0	&	0	&	0.91	&	0.11	&	0.25	&	0.13	&	0	&	1.	&	0.88	 \\
9	&	0	&	0	&	0	&	0	&	0.39	&	0.29	&	1.	&	0.78	&	0	&	0.62	 \\
10	&	0	&	0	&	0	&	0	&	0	&	0.41	&	1.	&	0.55	&	0.68	&	0	 \\  \hline
 \end{tabular}
 }
\vspace{0.5cm}
		\centering
\caption{\label{test1-res} Test 1}
{\small
\begin{tabular}{|c||c|c|c||c|c|   }
 \hline
method             &  \multicolumn{3}{c||}{MLSM}            &  \multicolumn{2}{c|}{  MLSM-AA(1)}  \\ \hline
$k_{max}$       &  \multicolumn{2}{c|}{1}     &  2        &    \multicolumn{2}{c|}{1}           \\ \hline
$s_{max}$       &  1        &  2                         & 1          & 1       &   2             \\ \hline
$N_t$               &   16     &  15                       & 15        &     15    &     15            \\
$\rho_{num}$  & 0.20      &  0.20                   &  0.19    &     0.20       &     0.20         \\
$M_{lo}$          & 2         &  3                         &  4          &   2       &   3          \\ \hline
 \end{tabular}
 }
\end{table}

\begin{figure}[t]
	\centering \hspace*{-.5cm}
	\subfloat[\label{test1-hist} Test 1]{\includegraphics[width=.5\textwidth]{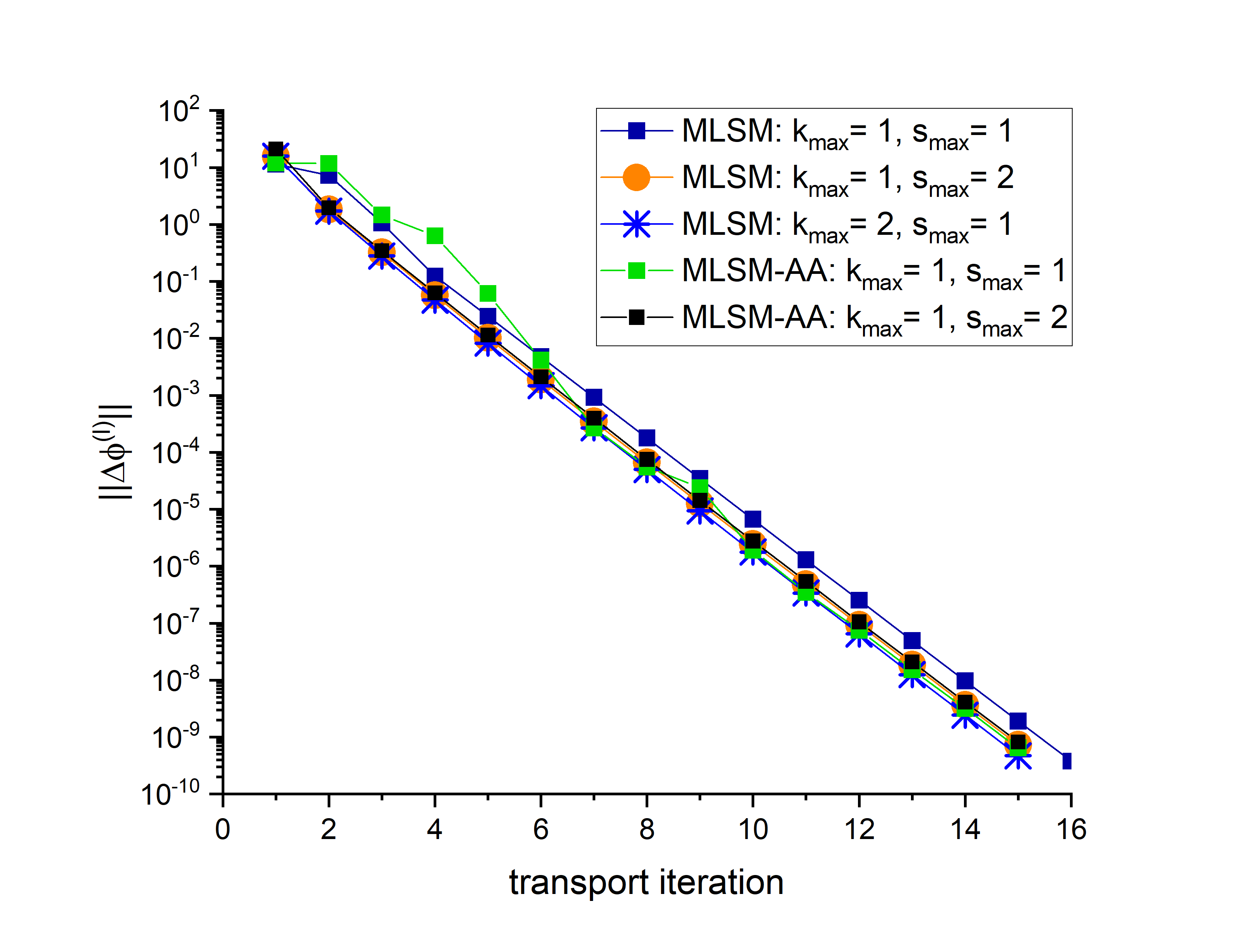}}
	\subfloat[\label{test2-mlsm-hist-1} Test 2, MLSM, $k_{max}=1$]{\includegraphics[width=.5\textwidth]{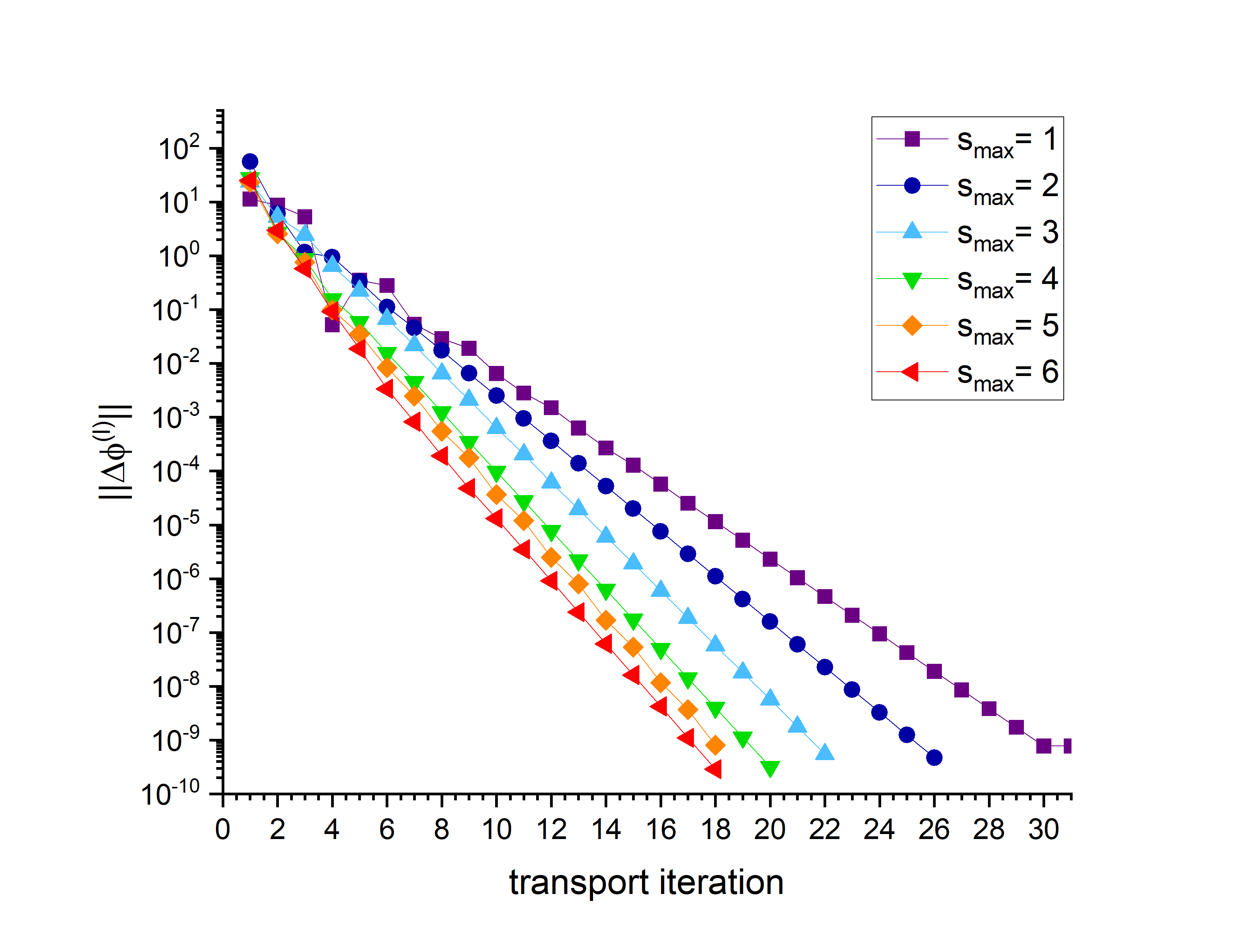}}\\
	\subfloat[\label{test2-mlsm-hist-2} Test 2, MLSM, $k_{max}=2,3,4,5$]{\includegraphics[width=.5\textwidth]{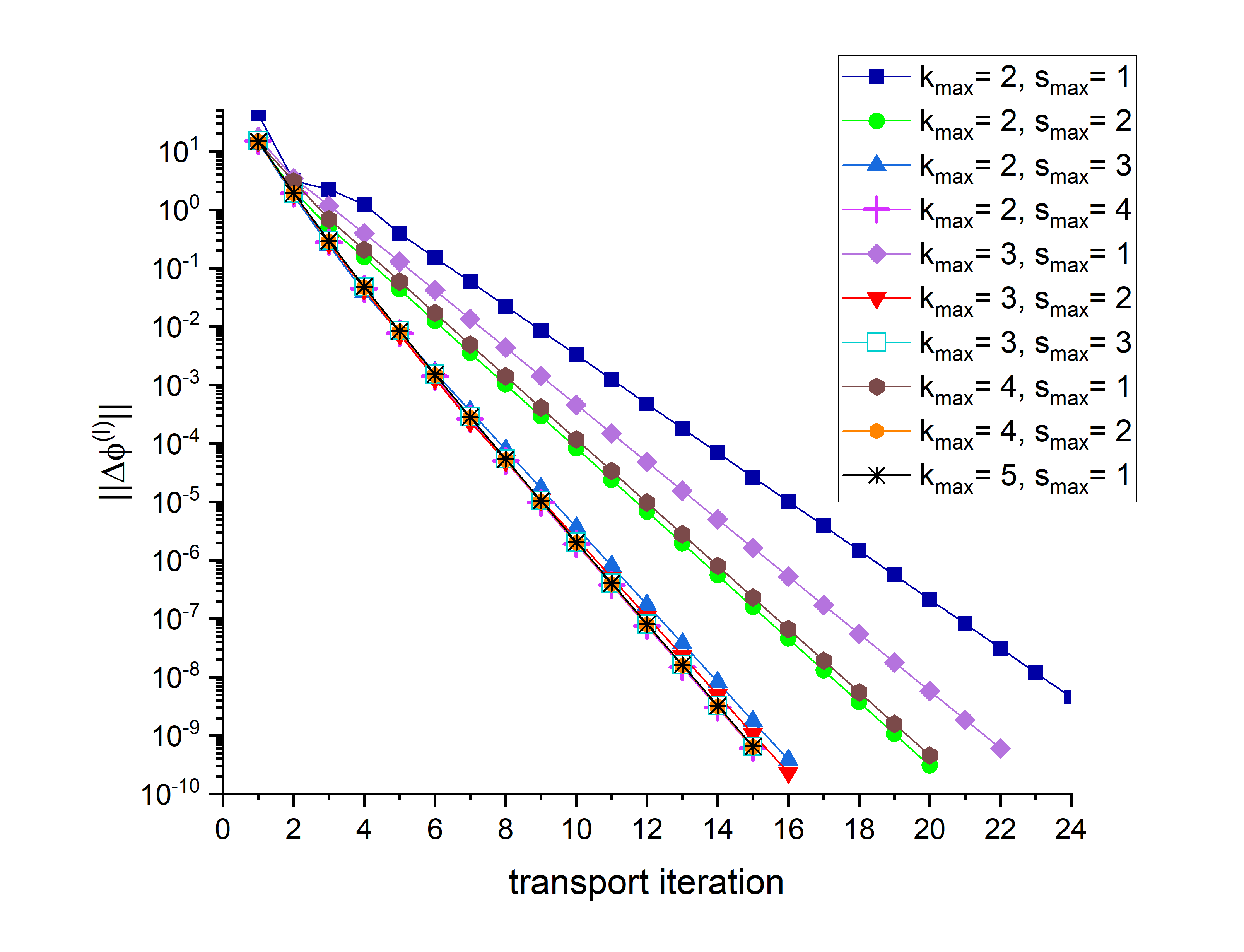}}
	\subfloat[\label{test2-mlsmaa-hist} Test 2,  MLSM-AA(1), $k_{max}=1,2,3$]{\includegraphics[width=.5\textwidth]{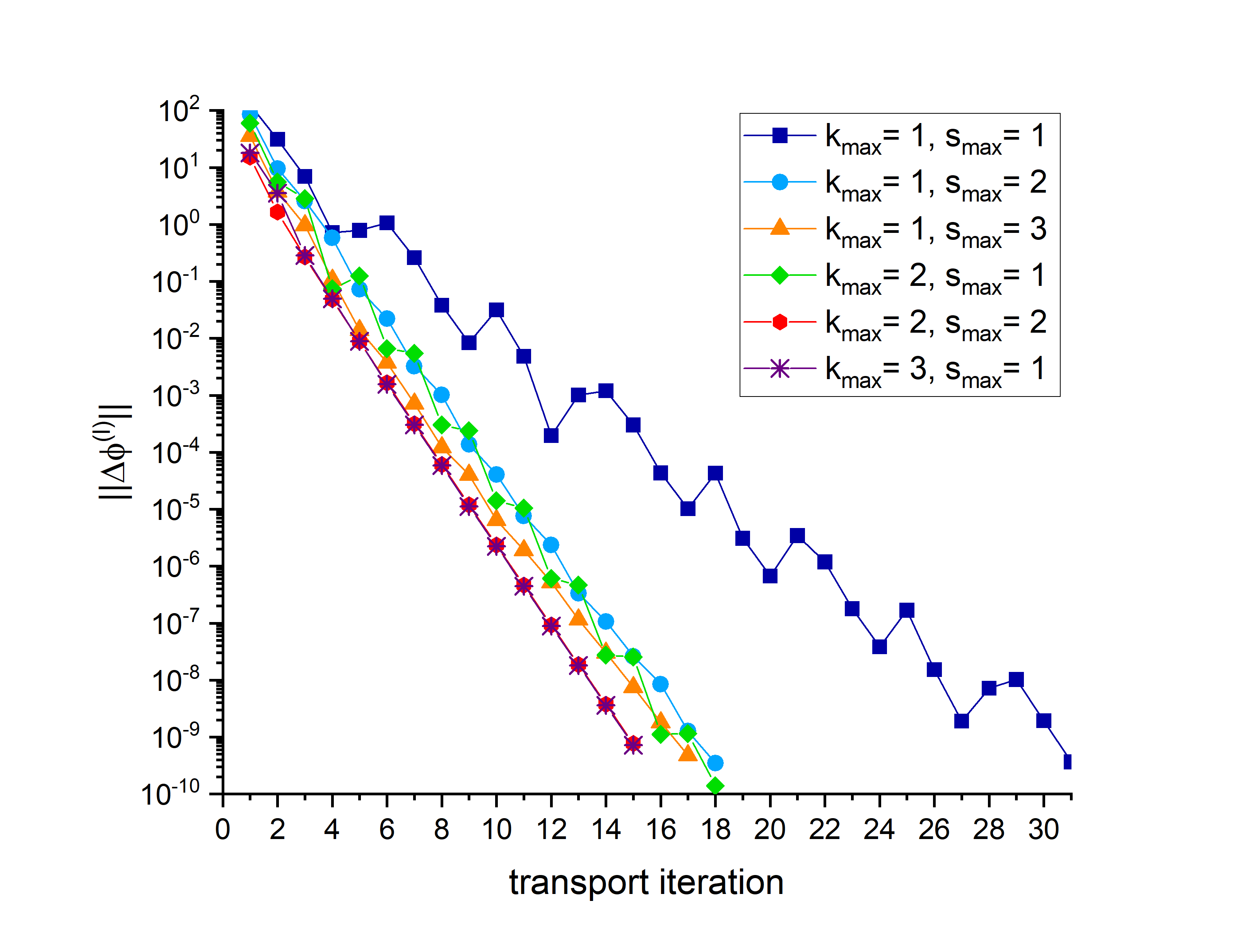}}
	\caption{\label{test2} Residual histories.}
\end{figure}

 Table  \ref{test1-res}  shows the numbers of outer transport iterations ($N_t$)
and numerically estimated spectral radii ($\rho_{num}$) for the MLSM and  MLSM-AA(1) methods
based on the rates of convergence during last iterations.
 The residual histories for the  methods are presented in Figures \ref{test1-hist}.
 The Fourier analysis in continuous form yields that the  value  of theoretical spectral radius (for an infinite-medium problem) of
 source iterations (SI) in this test is $\rho_{th}^{SI}=0.96$.
The full DSA  (FDSA) method has  $\rho_{th}^{FDSA}=0.21$ \cite{jmc-jsw-dya-jhc-m&c2021}.
The study of  the grey DSA (GDSA) and decoupled DSA (DDSA) showed that in this problem
$\rho_{num}^{GDSA}=0.55$ and $\rho_{num}^{DDSA}=0.26$ \cite{jmc-jsw-dya-jhc-m&c2021}.

On each   transport iteration, the MLSM algorithm executes
$k_{max}s_{max}$  parallel solves of   LOSM equations in groups and $k_{max}$ solves of grey LOSM equations. Thus, it performs  \linebreak $M_{lo}=k_{max}(s_{max} + 1)$ low-order solves
where each solve of group LOSM equations is accounted as one because of parallel execution of groups.
 This measure can be
used to evaluate the algorithm efficiency for the given number of transport iterations.

The results show that   MLSM with $k_{max}=1$ and $s_{max}=1$ converges
fast. Just one extra cycle over group LOSM equations ($s_{max}=2$)  leads to
 the target number of transport iterations ($N_t =15$).
 The estimated spectral radius of this algorithm  is $\rho_{num}=0.2$.
 This version of the algorithm has the smaller   number of cycles of low-order solves ($M_{lo}=3$)
  compared to the algorithm with $k_{max}=2$ and $s_{max}=1$.
The MLSM-AA(1) method  slightly  improves convergence in this test.
This algorithm with $k_{max}=1$ and $s_{max}=1$
shows the best performance in this test.
It converges in $N_t =15$  requiring   $M_{lo}=2$.

{\bf Test 2}.
This problem is similar to Test 1. It is defined with  the  moderator material from C5G7 benchmark with 7-group cross sections \cite{c5g7}.  Table \ref{test2-xs} shows the   cross section (see \ref{appendix}).
The  connection strength of groups  is presented in Table \ref{test2-strength}.
The groups are strongly connected to neighbouring groups.
The connection  to  distant groups   is  very weak.
The group scattering is  very  high in all groups ($0.985949 \le c_g \le  0.999961$).
The target number of outer transport iterations is equal to 15.

\begin{table}[t]
		\centering
\caption{\label{test2-strength}  Connection strength of  groups ($S_{gg'}$)  in Test 2}
{\small
\begin{tabular}{|c|c|c|c|c|c|c|c|}
\hline
$g \, \backslash \, g'$	& 1 & 2 & 3 & 4 & 5 & 6 & 7 \\ \hline
1	& 0 & 0 & 0  & 0  & 0  & 0  & 0  \\
2	& 1.   & 0 & 0  & 0  & 0  & 0  & 0  \\
3	& 5.5$\!\times\! 10^{-3}$ & 1. & 0 & 0 & 0 & 0 & 0 \\
4	& 1.7$\!\times\! 10^{-5}$& 2.8$\!\times\! 10^{-3}$ & 1. & 0 & 3.2$\!\times\! 10^{-4}$ & 0 & 0 \\
5	& 1.3$\!\times\! 10^{-7}$ & 1.2$\!\times\! 10^{-4}$ & 4.1$\!\times\! 10^{-2}$ & 1. & 0 & 5.3$\!\times\! 10^{-3}$ & 0 \\
6	& 0 & 1.5$\!\times\! 10^{-5}$ & 5.2$\!\times\! 10^{-3}$ & 1.3$\!\times\! 10^{-1}$ & 1. & 0 & 2.6$\!\times\! 10^{-1}$ \\
7	& 0 & 2.0$\!\times\! 10^{-6}$ & 9.4$\!\times\! 10^{-4}$ & 2.3$\!\times\! 10^{-2}$ & 1.1$\!\times\! 10^{-1}$ & 1. & 0 \\ \hline
 \end{tabular}
 }
\vspace{0.5cm}
		\centering
\caption{\label{test2-mlsm-res} Test 2:  MLSM}
\resizebox{1\textwidth}{!}{
\begin{tabular}{|c||c|c|c|c|c|c||c|c|c|c||c|c|c||c|c||c|   }
 \hline
 $k_{max}$     &  \multicolumn{6}{c||}{1}                        &   \multicolumn{4}{c||}{2}   &      \multicolumn{3}{c||}{3}    &      \multicolumn{2}{c||}{4}   &      5  \\ \hline
$s_{max}$      &  1    &  2     & 3      &  4     &  5     & 6     &   1     & 2      & 3       & 4     &  1    &  2      & 3          & 1      &  2       & 1 \\ \hline
$N_t$              &  31  &  26   & 22    &  20    &  18   & 18   &   26   & 20    & 16     & 15    & 22   & 16     & 15        &  20   &  15     & 15 \\
$\rho_{num}$ & 0.45 & 0.38 & 0.31 &  0.28 & 0.22 & 0.26 &  0.38 & 0.29 &  0.22 & 0.20 & 0.33 & 0.22 &  0.20    &  0.29 & 0.20   &  0.20\\
$M_{lo}$          &  2   & 3      & 4      &  5     & 6      & 7      &  4     & 6     &  8      & 10     &  6    &  9     &    12     &  8      & 12     &  10 \\  \hline
 \end{tabular}
 }
  \vspace{0.5cm}
		\centering
\caption{\label{test2-mlsmaa-res} Test 2:  MLSM-AA(1)}
{\small
\begin{tabular}{|c||c|c|c||c|c||c|    }
 \hline
  $k_{max}$     &  \multicolumn{3}{c||}{1}               &  \multicolumn{2}{c||}{2}   &  3   \\ \hline
$s_{max}$       &  1    & 2     & 3                                &  1   &  2                            &  1  \\ \hline
$N_t$              &  31   & 18   & 17                              &  18 & 15                           &   15    \\
$\rho_{num}$  &  n/a &  0.27 & 0.26                           & n/a &  0.20                       &  0.20  \\
$M_{lo}$          &  2   &   3    &   4                             &   4     & 6                            &   6   \\  \hline
 \end{tabular}
 }
\end{table}
The numbers of outer transport iterations  and numerically estimated spectral radii   for the MLSM and  MLSM-AA(1) methods are listed in Tables \ref{test2-mlsm-res} and \ref{test2-mlsmaa-res}, respectively.
The residual histories for both methods are presented in Figures  \ref{test2-mlsm-hist-1}-\ref{test2-mlsmaa-hist}.
The theoretical spectral radii of SI and FDSA   are
 $\rho_{th}^{SI}=0.98$  and  $\rho_{th}^{FDSA}=0.22$, respectively.
The analysis of   GDSA  and  DDSA  showed that in this problem
$\rho_{num}^{GDSA}=0.67$ and $\rho_{num}^{DDSA}=0.32$ \cite{jmc-jsw-dya-jhc-m&c2021}.

The results show that the MLSM algorithm with $k_{max}=2$ and $s_{max}=4$ converges in $N_t=15$ requiring $M_{lo}=10$  per transport iteration.  This method with $k_{max}=5$  with  only $s_{max}=1$  also converges in $N_t=15$ and needs the same number $M_{lo}$ per transport iteration.
Application of Anderson acceleration significantly affects permeance of  the MLSM method in this test.
The most efficient is  the MLSM-AA(1) algorithm with $k_{max}=2$ and $s_{max}=3$
 that executes $M_{lo}=6$ per transport iteration.
This algorithm converges steadily with  estimated spectral radius   $\rho_{num}=0.20$.
We note  that the MLSM-AA(1) method with $k_{max}=1$ and $s_{max}=1$
  showed irregular convergence behavior. This is the effect of using in Anderson acceleration
   just one residual    of the  solution of the group LOSM equations  for $s=1$  and
    the residual of the initial guess ($s=0$) that is the high-order solution from the transport sweep.
The trace of this effect can be also noticed in  convergence behaviour of the MLSM-AA(1) method with $k_{max}=2$ and $s_{max}=1$.

\section{\label{sec:conc} Conclusions}

We developed new multilevel iterative methods for fixed-source multigroup particle transport problems
that can be applied for parallel computations.
Numerical results are promising. They show that the algorithms accelerate iterative convergence
and effectively solve multigroup test problems with  down- and upscattering as well as with high scattering ratios in groups.
 More analysis is needed to study  properties of MLSM iterative algorithms.
 Further work will include extension to multi-D geometries and application of more general version of Anderson acceleration.
This kind of transport algorithms for parallel computations can also be developed on the basis of the quasidiffusion (VEF) method \cite{dya-vyag-nse-2011,dya-jsw-nse-2018}.

\section*{Acknowledgements}

Los Alamos Report LA-UR-20-26669.
This work was funded in part by Los Alamos National Laboratory, which is operated by Triad National Security, LLC,
for the National Nuclear Security Administration of U.S. Department of Energy (Contract No. \linebreak 89233218NCA000001).
The work of the second author (JMC) was funded by LANL through a summer research internship in the CCS-2 group.

\bibliographystyle{elsarticle-num}
\bibliography{MLSM-mc2021}

\clearpage
\appendix
\gdef\thesection{Appendix \Alph{section}}
\section{Cross Section Data}
\label{appendix}
\begin{table}[h!]
		\centering
\caption{\label{test1-xs} Cross section data for Test 1 \cite{jmc-jsw-dya-jhc-m&c2021}}
\resizebox{1\textwidth}{!}{
\begin{tabular}{|c|c|c|c|c|c|c|c|c|c|c|}
\hline
 $ g$ &  1 &  2 & 3  & 4 & 5 &  6 &  7 &  8 &  9 &  10  \\ \hline
 $\sigma_{t,g}$ &  2.49756 &  2.01650 & 1.51992  & 1.67388 & 2.36661 &  1.50008 &  2.37543 &  2.36241 &  2.04640 &  1.59740  \\
 $c_{g}$ &    0.979581 &   0.944816 &   0.952295 &   0.926035 &   0.978471 &   0.9  &   0.987210 &   0.9999  &   0.904252 &   0.966192  \\
\hline \hline
$\tilde g$    &   $\sigma_{s,1 \to \tilde g}$  &      $\sigma_{s,2 \to \tilde g}$  &  $\sigma_{s,3 \to \tilde g}$  &  $\sigma_{s,4 \to \tilde g}$  & $\sigma_{s,5 \to \tilde g}$  & $\sigma_{s,6 \to \tilde g}$  & $\sigma_{s,7 \to \tilde g}$  &      $\sigma_{s,8 \to \tilde g}$  &     $\sigma_{s,9 \to \tilde g}$  &   $\sigma_{s,10 \to \tilde g}$     \\ \hline
1  &     0.835282  &   0  &   0  &   0  &   0  &   0  &   0  &   0  &   0  &   0       \\
2  &     0.401686  &    0.566521 &     0  &   0  &   0  &   0  &   0  &   0  &   0  &   0      \\
3  &     0.404298  &    0.569454 &   0.420634  &   0   &   0   &   0   &   0   &   0    &   0    &   0    \\
4  &     0.498922  &    0.264139 &   0.179242 &     0.0828011     &   0   &   0    &   0    &   0    &   0    &   0    \\
5  &     0.306376  &    0.0657747 &    0.148397  &   0.307318  &   1.30088    &   0   &   0    &   0    &   0    &   0     \\
6  &      0.439338  &    0.362807 &    0.564376 &    0.456018  &   0.0715262   &   0    &   0    &   0    &   0  & 0    \\
7  &                   0 &     0  &    0.336331  &  0.122044  &  0.259295  &   0.623241 &   0.812409 &   1.28728  &   0.278371 &  0.301517    \\
8  &       0 &    0 &    0 &    0.473528  &   0.0566290  &    0.128925  &    0.0676741  &   0.123057  &    0.518149  &   0.457140     \\
9  &       0    &   0   &   0    & 0      & 0.242843  &    0.180473  &   0.622078  &    0.485474  &    0.483321  &    0.386770     \\
10  &     0    &   0    &   0    &   0   &   0  &  0.345904  &    0.842890  &    0.466367  &    0.570623  &    0.397965     \\  \hline
\end{tabular}
}
\vspace{0.5cm}
 		\centering
\caption{\label{test2-xs} Cross section data for Test 2 \cite{c5g7}}
\resizebox{1\textwidth}{!}{
\begin{tabular}{|c|c|c|c|c|c|c|c|}
\hline
 $ g$  &  1 &  2 & 3  & 4 & 5 &  6 &  7    \\ \hline
 $\sigma_{t,g}$ & 0.159206 & 0.412970 & 0.590310 & 0.584350 & 0.718000 &  1.25445 & 2.65038 \\
 $c_g$ &   0.996225 & 0.999961 & 0.999429 & 0.996679 & 0.992003 & 0.988042 & 0.985949   \\ \hline
\hline
$\tilde g$    &   $\sigma_{s,1 \to \tilde g}$  &      $\sigma_{s,2 \to \tilde g}$  &  $\sigma_{s,3 \to \tilde g}$  &  $\sigma_{s,4 \to \tilde g}$  & $\sigma_{s,5 \to \tilde g}$  & $\sigma_{s,6 \to \tilde g}$  & $\sigma_{s,7 \to \tilde g}$  \\ \hline
1 &    4.44777$\!\times \! 10^{-2}$ &    0  &     0  &      0  &      0  &      0  &      0     \\
2 &      1.134$\!\times \! 10^{-1}$ &    2.82334$\!\times \! 10^{-1}$ &    0  &     0  &      0  &     0  &      0     \\
3 &      7.2347$\!\times \! 10^{-4}$ &    1.2994$\!\times \! 10^{-1}$ &    3.45256$\!\times \! 10^{-1}$ &    0  &      0  &      0  &      0     \\
4 &       3.7499$\!\times \! 10^{-6}$ &    6.234$\!\times \! 10^{-4}$ &    2.2457$\!\times \! 10^{-1}$ &    9.10284$\!\times \! 10^{-2}$ &    7.1437$\!\times \! 10^{-5}$ &    0  &      0     \\
5 &      5.3184$\!\times \! 10^{-8}$ &    4.8002$\!\times \! 10^{-5}$ &    1.6999$\!\times \! 10^{-2}$ &    4.1551$\!\times \! 10^{-1}$ &   1.39138$\!\times \! 10^{-1}$ &    2.2157$\!\times \! 10^{-3}$ &    0     \\
6 &       0  &      7.4486$\!\times \! 10^{-6}$ &    2.6443$\!\times \! 10^{-3}$ &    6.3732$\!\times \! 10^{-2}$ &    5.1182$\!\times \! 10^{-1}$ &    6.99913$\!\times \! 10^{-1}$ &    1.3244$\!\times \! 10^{-1}$      \\
7 &       0  &      1.0455$\!\times \! 10^{-6}$ &    5.0344$\!\times \! 10^{-4}$ &    1.2139$\!\times \! 10^{-2}$ &    6.1229$\!\times \! 10^{-2}$ &    5.3732$\!\times \! 10^{-1}$ &    2.4807    \\ \hline
\end{tabular}
}
\end{table}

\end{document}